\documentclass[11pt]{article}
\usepackage{amsmath,amsthm,amsfonts}

\newcommand{\LL}{\mathfrak{L}}
\newcommand{\RR}{\mathfrak{R}}
\newcommand{\MM}{\mathcal{M}}
\renewcommand{\=}{\doteq}

\newtheorem{thm}{Theorem}[section]
 \newtheorem{cor}[thm]{Corollary}
 \newtheorem{prop}[thm]{Proposition}

\theoremstyle{definition}

\theoremstyle{definition}
 \newtheorem{rem}[thm]{Remark}
\numberwithin{equation}{section}
\begin{document} 
\title{\bf Moufang symmetry V.\\Triple closure}
\author{Eugen Paal}
\date{}
\maketitle

\thispagestyle{empty}

\begin{abstract}
Triple closure of the infinitesimal translations of an analytic Moufang loop is inquired. This property is equivalent to reductivity and relates Mal'tsev algebras to the Lie triple systems.
\end{abstract}

\section{Introduction}

In the present paper we inquire triple closure of the infinitesimal translations of a local analytic Moufang loop. This property is equivalent to reductivity and relates Mal'tsev algebras to the Lie triple systems. The paper can be seen as a continuation of \cite{Paal1,Paal2,Paal3,Paal4}.

\section{Reductivity and triple closure}

We know from \cite{Paal4} that the infinitesimal translations $L_x$, $R_x$ and $M_x$ ($x\in T_e(G)$) of  a local analytic Moufang loop $G$ satisfy the \emph{reductivity} conditions
\begin{subequations}
\label{red_lrm}
\begin{align} 
[Y(x;y),L_z]&=L_{[x,y,z]}\\
[Y(x;y),R_z]&=R_{[x,y,z]}\\
[Y(x;y),R_z]&=R_{[x,y,z]}
\end{align}
\end{subequations}
where the Yamagutian $Y$ and the \emph{Yamaguti brackets} $[\cdot,\cdot,\cdot]$ are given by
\begin{align*}
6Y(x;y)&=[L_x,L_y]+[R_x,R_y]+[M_x,M_y]\\
[x,y,z]
&=6(x,y,z)+2[[x,y],z]\\
&=[x,[y,z]]-[y,[x,z]]+[[x,y],z]
\end{align*}
Now define  \cite{Loos} the \emph{Loos brackets}  $\{\cdot,\cdot,\cdot\}$ by
\begin{subequations}
\label{loos}
\begin{align} 
3\{x,y,z\}
&\=6(x,y,z)+3[[x,y],z]\\
&=[x,y,z]+[[x,y],z]\\
&=[x,[y,z]]-[y,[x,z]]+2[[x,y],z]
\end{align}
\end{subequations}
Since $(L,R)$ is a Moufang-Mal'tsev pair we know from \cite{Paal2}

\begin{prop} 
Let $(S,T)$ be a Moufang-Mal'tsev pair. Then 
\begin{subequations}
\label{y-lrm}
\begin{align}
6[Y(x;y),L_{z}]&=3[[L_{x},L_{y}],L_{z}]-L_{[[x,y],z]}\\
6[Y(x;y),R_{z}]&=3[[R_{x},R_{y}],T_{z}]-R_{[[x,y],z]}\\
6[Y(x;y),M_{z}]&=3[[M_{x},M_{y}],M_{z}]-M_{[[x,y],z]}
\end{align}
\end{subequations}
for all $x,y,z$ in $M$.
\end{prop}

\begin{thm}[triple closure]
The infinitesimal translations $L_x$, $R_x$ and $M_x$ ($x\in T_e(G)$) of  a local analytic Moufang loop $G$ satisfy the \emph{triple closure} conditions 
\begin{subequations}
\label{triple_lrm}
\begin{align} 
[L_x,L_y],L_z]&=L_{\{x,y,z\}}\\
[R_x,R_y],R_z]&=R_{\{x,y,z\}}\\
[M_x,M_y],M_z]&=M_{\{x,y,z\}}
\end{align}
\end{subequations}
These relations hold if and only if  the reductivity conditions hold.
\end{thm}

\begin{proof}
Use (\ref{y-lrm}a-c) and reductivity (\ref{red_lrm}a-c).
\end{proof}

The triple close property (\ref{triple_lrm}a-c) means that the vector spaces 
\begin{equation*}
\LL \=\{L_x|x\in T_e(G)\}, \quad \RR \=\{R_x|x\in T_e(G)\},\quad \MM \=\{M_x|x\in T_e(G)\}
\end{equation*}
are closed under the double Lie bracketing. N.~Jacobson \cite{Jac} called such spaces (subspaces of the Lie algebras) the \emph{Lie triple systems}. Thus spaces $L$, $R$ and $M$ turn out to be the Lie triple systems.

\begin{rem}
The Lie triple systems turn out to be \cite{Yam57} the tangent algebras of the symmetric spaces. In \cite{SabMih} it was shown that an analytic Moufang loop is an affine symmetric space such that its motion group is generated by the left and right translations of the given analytic Moufang loop. Thus it is not surprising that the spaces  $L$, $R$ and $M$ are the Lie triple systems.
\end{rem}

\begin{thm}
The infinitesimal translations of an analytic Moufang loop $G$ satisfy the ralations
\begin{subequations}
\label{sym_lrm}
\begin{align} 
L(x;y)&\=[L_x,L_y]\\
[L(x;y),L_z]&=L_{\{x,y,z\}}\\
[L(x;y),L(z;w)]&=L(\{x,y,z\};w)+L(z;\{x,y,w\}
\end{align}
\end{subequations}
\end{thm}

\begin{proof} 
Denote $L(x;y)\=[L_x,L_y]$ and consider the Jacobi identity
\begin{align*}
0
&=[[[L_x,L_y],L_z],L_w]+[[L_z,L_w],[L_x,L_y]]+[[L_w,[L_x,L_y]],L_z)\\
&=[L_{\{x,y,z\}},L_w]-[L(x;y),L(z;w)]+[L_z,L_{\{x,y,w\}]}]
\tag*{\qed}
\end{align*}
\renewcommand{\qed}{}
\end{proof}

We can see that (\ref{sym_lrm}a-c) are the structure relations for the symmetric space $\LL+[\LL,\LL]$. By triality, the similar structure relations hold for $\RR+[\RR,\RR]$ and $\MM+[\MM,\MM]$.

It is natural that triple closure property of the spaces $\LL$, $\RR$ and $\MM$ reflects in the corresponding properties of the Loos bracket.
\begin{thm}
The Loos brackets $\{\cdot,\cdot,\cdot\}$ satisfy relations
\begin{subequations}
\label{loos_properties}
\begin{gather} 
\{x,y,z\}=-\{y,x,z\}\\
\{x,y,z\}+\{y,z,x\}+\{z,x,y\}=0\\
\{x,y,\{z,w,v\},\}=\{\{x,y,z\},w,v\}+\{z,\{x,y,w\},v\}+\{z,w,\{x,y,v\}\}
\end{gather}
\end{subequations}
\end{thm}

\begin{proof}
Relation (\ref{loos_properties}a) is evident from (\ref{loos}a-c). Relation (\ref{loos_properties}b) follows easily by using the triple closure from the Jacobi identity
\begin{equation*}
[[L_x,L_y],L_z]+[[L_y,L_z],L_x]+[[L_z,L_x],L_y]=0
\end{equation*}
Finally, (\ref{loos_properties}c) follows by using the triple close from the Jacobi identity
\begin{equation*}
[[[L_x,L_y],[L_z,L_w]],L_v]+[[[L_z,L_w],L_v],[L_x,L_y]]+[[L_v,[L_x,L_y]],[L_z,L_w]]=0
\tag*{\qed}
\end{equation*}
\renewcommand{\qed}{}
\end{proof}

\begin{cor}[see also \cite{Loos}]
The tangent algebra $\Gamma$ of a local analytic Moufang loop is  a Lie triple system with the Loos brackets (\ref{loos}a-c).
\end{cor}

\section*{Acknowledgement}

Research was in part supported by the Estonian Science Foundation, Grant 6912.

\bigskip\noindent
Department of Mathematics\\
Tallinn University of Technology\\
Ehitajate tee  5, 19086 Tallinn, Estonia\\
E-mail: eugen.paal@ttu.ee

\end{document}